\documentclass[12pt,reqno]{article}

\usepackage{amsmath}
\usepackage{amsthm}
\usepackage{amscd}
\usepackage{amssymb}
\usepackage{graphicx}
\usepackage{epsf}
\theoremstyle{plain}
\newtheorem{theorem}{Theorem}[section]
\newtheorem{proposition}[theorem]{Proposition}

\newtheorem{corollary}[theorem]{Corollary}

\newtheorem{conjecture}[theorem]{Conjecture}
\newtheorem{question}[theorem]{Question}

\newtheorem{definition}[theorem]{Definition}
\newtheorem{remark}[theorem]{Remark}

\def\int{\text{int}}
\def\invlimit{\smash{\lim\limits_{\raise1pt\hbox{$\longleftarrow$}}}\vphantom{\big(}}
\def\inter{\hskip 1.5pt\raise4pt\hbox{$^\circ$}\kern -1.6ex}
\def\skel(#1,#2){#1^{(#2)}}
\def\hyp {\hbox {\rm {H \kern -2.8ex I}\kern 1.25ex}}
\def\reals {\hbox {\rm {R \kern -2.8ex I}\kern 1.15ex}}
\def\integers {\hbox {\rm { Z \kern -2.8ex Z}\kern 1.15ex}}
\def\naturals {\hbox {\rm {N \kern -2.8ex I}\kern 1.20ex}}
\def\rationals {\hbox {\rm { Q \kern -2.2ex l}\kern 1.15ex}}
\def\hyp {\hbox {\rm {H \kern -2.7ex I}\kern 1.25ex}}

\begin{document}

\title{Connected sums of knots and weakly reducible Heegaard splittings}

\author{Yoav Moriah \thanks{Supported by The Fund for Promoting Research at the Technion 
grant 100-053}}

\date{}

\maketitle

\begin{abstract}
This paper studies the question of whether minimal genus Heegaard splittings of exterior spaces of 
knots which are connected sums are weakly reducible or not. Furthermore it is shown that the Heegaard 
splittings of the knots used by Morimoto to show that tunnel number can be sub-additive are all strongly 
irreducible. These are the first  examples of strongly irreducible minimal genus Heegaard splittings of 
composite knots. We also give a characterization of when is a set of primitive annuli on a handlebody
simultaneously primitive. This  characterization  is different from that given in [Go]. 
\end{abstract}

\vskip20pt

\section{Introduction}
\label{intro}

\vskip20pt

For some time it is known that there is a connection between the existence of  closed incompressible 
surfaces in a $3$-manifold and the nature of its Heegaard splittings.  See for example [CG], 
 [Ha], [LM], [Mh], [Mo3], [Mo4], [MS], [Sc]. In this paper we begin to explore this 
connection with respect to essential surfaces with boundary and as a first step study spaces containing
essential annuli. A special case of manifolds which contain essential (i.e., incompressible non-boundary 
parallel) annuli are exterior spaces of connected sums of  knots in $S^3$. These manifolds are obtained 
from two knot exterior spaces by gluing them together along a meridional annulus $A$.

Given a Heegaard  splitting for a manifold $M$ (i.e., a decomposition $M = V_1 \cup V_2$, $ V_1 \cap V_2 = \Sigma$, 
where $V_i$, $ i = 1,2$ are  compression bodies  and  $\Sigma = \partial V_1 = \partial V_2$ is the Heegaard  
surface)  let ${\cal C}_i$ denote  the set of all essential simple curves on $\Sigma$ which  bound disks 
in $V_i$. Define: $d(V_1, V_2) = min\{d(C_1, C_2)| C_i \in{ \cal C}_i(\Sigma)\},$ where $d(C_1, C_2)$  is 
measured in the curve complex $\cal C$ of $\Sigma$. In particular a Heegaard splitting will be reducible if 
$d(V_1,V_2) = 0$, weakly reducible if $d(V_1, V_2) \leq 1$ and strongly irreducible if $d(V_1,V_2) \geq 2$.
Note that any knot exterior of a knot which is a connected sum, contains at least two essential tori. It is a result 
of Hempel [He] and Thompson  [Th] that if a manifold contains an essential torus then any Heegaard splitting 
$(V_1, V_2)$ has $d(V_1, V_2) \leq2$. In general we have a result by Hartshorn [Ha] that if an irreducible 
$3$-manifold $ M$ contains a closed incompressible surface of genus $g$ the distance  of any  Heegaard splitting 
$(V_1, V_2)$  of $M$ is less than or equal to $2g$. 

As the Euler characteristic of  an annulus is $0$, just like that of a torus, and it is also a twice punctured 
$2$-sphere one might \lq\lq hope" that the theorem of Hartshorn might be extended to say that an 
essential annulus in a $3$-manifold with torus boundary will imply that $d(V_1, V_2) \leq 1$. In other 
words: Any Heegaard splitting of such a manifold will be weakly reducible. Evidence in this direction 
is in  [LM] where the authors describe a very large class  of knots in $S^3$ for which the connected 
sum yields manifolds with a minimal genus Heegaard  splitting which are weakly reducible. In this direction 
we prove the following theorems:

\vskip10pt

\noindent 
{\bf Theorem \ref{Supaddthm}} {\it Given knots $K_1, K_2$  and  $ K = K_1 \# K_2$ in $S^3$ 
for which the tunnel number satisfies $t(K) = t(K_1) + t(K_2) +1$ i.e., $t(K)$ is super additive, 
then there is a minimal genus Heegaard splitting of  $E(K)$ which is  weakly reducible.}

\vskip10pt

\noindent {\bf Theorem \ref{wredProp}} {\it Let $K_1, K_2$  and  $ K = K_1 \# K_2$ be knots  in $S^3$
and $(V_1^i , V_2^i ), i =1,2$ be  Heegaard splittings  for $E(K_i)$. If $(V_1^1, V_2^1)$  and $(V_1^2 , V_2^2)$
induce a  Heegaard splitting $(V_1, V_2)$  of  $E(K)$  then $(V_1, V_2)$ is a weakly reducible Heegaard splitting.}

\vskip12pt

In particular this theorem says that if one of $E(K_1)$ or $E(K_2)$ has a $\mu$-primitive minimal genus
Heegaard splitting then $E(K_1 \# K_2)$ will have a weakly reducible Heegaard splitting of genus $g_1 + g_2 - 1$,
where $g_i = genus(E(K_i)) $

\vskip12pt

\noindent {\bf Theorem \ref{nhorwThm}}  {\it Let $K = K_1 \# K_2 \subset S^3$ be a  knot. 
Any Heegaard surface $\Sigma$ for $E(K)$ which does not contain any $\Sigma$ horizontal surfaces 
is weakly reducible.}

\vskip10pt

\noindent Finally:

\vskip10pt

\noindent {\bf Theorem \ref{singtunthm}}  {\it  Let $K_1, K_2$ be prime knots in $S^3$ and  
$K = K_1 \# K_2$. Assume that $t(K) = t(K_1) + t(K_2)$ and  $t(K_i) \leq 2 $.  Furthermore, 
assume that a minimal tunnel system for $K$ minimaly intersects a   decomposing annulus $A$ 
in a single point, then there is a Heegaard splitting of $E(K)$  of minimal genus  which is weakly reducible.}

\vskip10pt

However the connection between the distance of Heegaard splittings and the existence of an essential 
annulus is more complicated as shown by the following theorem. Let $K_n$ denote the knots as in [Mo2] : 

\vskip10pt

\noindent  {\bf Theorem 5.1.}  {\it Let  $K_n \subset S^3$ be the knot as in Fig.6 and
$K({ \frac\alpha  \beta})  \subset S^3$ a  2-bridge knot determined by 
$ { \frac\alpha  \beta} \subset \rationals $. Let $K$ denote the  connected sum $K_n  \#  K({ \frac\alpha  \beta})$, 
then the Heegaard splitting of  $E(K)$  determined by the minimal tunnel system for $K$, (as in Fig.6) is strongly
irreducible.}

\vskip10pt

These are the first examples of strongly irreducible Heegaard splittings of exteriors of connected sums.
These knots have the property that $g(E(K_1\#K_2))$ = $g(E(K_1)) + g(E(K_2)) - 2$, where g( ) denotes
the genus of the manifold in brackets. Hence a minimal genus Heegaard splitting of $E(K_1\#K_2)$
cannot possibly be induced by Heegaard splittings of the two knot spaces.

\vskip10pt

\noindent In light of the above I would like to propose the following conjecture:

\vskip10pt

\begin{conjecture}
\label{fstconj}
Given two knots $K_1, K_2$  in $S^3$ for which the tunnel number $t(K)$ satisfies 
$t(K_1 \# K_2) = t(K_1) + t(K_2)$,  then there is a minimal genus Heegaard splitting of  $E(K)$  
which is weakly reducible.

\end{conjecture}

\vskip10pt

The situation is further complicated by the possibility of a positive answer to the following open
question:

\vskip10pt

\begin{question}
Can a $3$-manifold $M$ have both weakly reducible and strongly irreducible
minimal genus Heegaard splittings ?
\label{qstn}
\end{question}

\vskip3pt

\noindent  For definitions of the above terminology see Sections \ref{prelim} and \ref{weak I}.

\vskip0pt

\noindent {\bf Acknowledgments:}  I would like to thank Ying-Qing Wu for suggesting the current
proof of Theorem \ref{simprim}. and other helpful remarks.

\vskip25pt

\section{Preliminaries}
\label{prelim}

\vskip15pt

\noindent 

Throughout the paper  $K_1$ and $K_2$ will be knots in $S^3$ and  $K = K_1 \# K_2$
will denote  the connected sum of $K_1$ and $K_2$.  The knots $K_i$ will be called  the 
{\it summands} of the {\it composite knot} $K$. Let $N()$ denote an open  regular 
neighborhood in $S^3$. An incompressible surface in a knot complement $E(K), K \subset S^3$ 
is called {\it meridional} if it has boundary  components which are meridian  curves of $\partial E(K)$. 

Recall that $(S^3,K)$ is obtained by removing  from each space  $(S^3,K_i),  i = {1,2},$ a  
small $3$-ball intersecting $K_i$ in a short unknotted arc and gluing the two remaining $3$-balls 
along  the $2$-sphere boundary so that the pair of points of $K_1$ on the $2$-sphere are identified 
with  the  pair of points of $K_2$. If we denote $S^3 - N(K)$ by $E(K)$ then $E(K)$ is obtained from 
$E(K_i), i = {1,2},$ by identifying a meridional annulus $A_1$ on $\partial E(K_1)$ with a 
meridional annulus $A_2$ on $\partial E(K_2)$. A knot $K \subset S^3$ is {\it prime} if it is not a 
connected sum of two non-trivial knots. The annulus $A_1 = A_2$ will be denoted by $A$ 
and called the{ \it decomposing annulus}. If both knots $K_1, K_2$ are prime then the decomposing 
annulus is unique up to isotopy  

 A {\it tunnel system} for an arbitrary knot $K \subset S^3$ is a collection of 
properly  embedded arcs  $\{t_1, \dots, t_n\}$  in  $S^3 - N(K)$  so that 
$S^3 - N(K \cup t_1 \cup \dots \cup t_n)$ is a handlebody.

Given a tunnel system for a knot $K \subset S^3$ note that the closure of 
$N(K \cup t_1 \cup \dots \cup t_n)$ is always a handlebody denoted by $V_1$ 
and the handlebody $S^3 - N(K \cup t_1 \cup  \dots \cup t_n)$ will be denoted by $V_2$. 
For a given knot $K \subset S^3$ the smallest cardinality of any tunnel system is called the 
{\it tunnel number} of $K$ and
is denoted by $t(K)$.

A compression body  $V$  is a  compact orientable and connected 3-manifold with a preferred 
boundary component $\partial_+V$  and is obtained from a collar of $\partial_+ V$ by attaching 
2-handles and 3-handles, so that the connected components of  $\partial_- V$ = $\partial V - \partial_+ V$ 
are all distinct from  $S^2$.  The extreme cases, where  $V$  is a handlebody i.e., $\partial_- V = \emptyset$,
or where $V = \partial_+V \times I$, are allowed.  Alternatively we can think of $V$ as obtained from
$(\partial_-V) \times I$ by attaching $1$-handles to $(\partial_-V) \times \{1\}$. An annulus in a 
compression body will be called a {\it spanning (or vertical) annulus} if it has one boundary component
on $\partial_+V$ and the other on $\partial_-V$.

Given a knot  $K \subset S^3$ a {\it Heegaard splitting } for $E(K)$ is a decomposition of $E(K)$ into
a compression body $V_1$ and a handlebody $ V_2 = S^3 - int(V_1)$. Hence, a tunnel system 
$\{t_1, \dots, t_n\}$ in $S^3 - N(K)$  for $K$ determines a Heegaard splitting of genus $n +1$ for $E(K)$.

When considering knot complements the operation of connected sum is well defined and not dependent on
the choice of the removed trivial ball pair $(B, t)$ as any two such ball pairs are isotopic in $E(K)$. However
when we are studying the additional structure of Heegaard splittings of composite knot complements we must be
careful as it is not clear that an isotopy of the ball pairs can induce an isotopy of the meridional annulus
preserving  the Heegaard surface. 
  
Given a Heegaard splitting $(V_1, V_2)$  for $ S^3 - N(K_1 \# K_2)$ we will choose  a 
decomposing annulus $A$  which intersects the compression body $V_1$ in two spanning annuli 
$A^*_1, A^*_2$  and a {\em minimal} collection of disks ${\cal D} = \{D_1, \dots, D_l\}$. Note also that 
$A$ intersects $V_2$ in a connected incompressible planar surface.

Let $ {\cal E} = \{E_1 , \dots, E_{t(K) + 1}\}$ be a complete meridian disk 
system for $V_2$, chosen to minimize the intersection ${\cal E}\cap A$. Since $V_2$ is a handlebody 
it is irreducible  and we can assume that no component of  ${\cal E} \cap A$ is a simple closed curve.

When we cut $E(K)$ along a decomposing annulus $A$ any Heegaard splitting $(V_1, V_2)$ of $E(K)$
induces  Heegaard splittings on both of $E(K_1)$ and $E(K_2)$, as follows: Set $V_1^i = (V_1 \cap E(K_i))
\cup_{{\cal D} \cup A_1^* \cup A_2^*} N(A)$, it is a compression body as it is a union of an $annulus \times I$
and some 1-handles along the two vertical annuli and a collection of disks. Now set  $V_2^i = V_2 - N(A)$, 
it is a handlebody since the annulus $A$  meets $V_2$ in an incompressible
connected planar surface $P$ which separates $V_2$ into two components each of which is a handlebody. 
Hence the pair $(V_1^i,V_2^i)$ is a Heegaard splitting for $E(K_i)$ and will be referred to as the {\it
induced Heegaard splitting} of $E(K_i)$.

We say that a curve on a handlebody is {\it primitive} if there is an essential disk in the handlebody 
intersecting the curve in a single point. An annulus $A$ on $H$ is primitive if its core curve is primitive.
A Heegaard splitting $(V_1, V_2)$ for $S^3 - N(K)$ will be called {\it $\mu$-primitive} if there is a 
spanning annulus $A\subset V_1$ such that $\partial A = \mu \cup \alpha$ where $\mu$ is a meridian 
and $\alpha$ is a primitive curve on $\partial V_2$.  Note that a curve on a handlebody $H$ is 
{\it primitive} if it represents a primitive element  in the free group $\pi_1(H)$.

Two Heegaard splittings $(V_1^i, V_2^i)$ for  $E(K_i)$ respectively, induce a decomposition of
$E(K)$ into $(V_1, V_2)$. We can think of $V_1^i$ as a union of $(\partial E(K_i) \times I)  \cup 
1-handles$, hence if we consider the ball pair $(B_i, N(t_i))$ and remove it from $E(K_i)$ we can think
of the decomposing annulus $A_i = \partial B_i - N(\partial t_i)$ as the union of two vertical annuli 
$A_1^{*i}, A_2^{*i}$ and a meridional annulus $A_i\subset \partial E(K_i) \times \{1\}  \subset 
\partial V_1^i = \partial V_2^i$. We obtain $V_1$ by gluing the compression bodies $V_1^1$ and 
$V_1^2$ along the two vertical annuli and $V_2$ by gluing $V_2^1$ and $V_2^2$ along a meridional 
annulus. Hence $V_1$ is always a compression body but $V_2$ is a handlebody if and only if the 
meridional annulus is a primitive annulus in $V_2^i$ for one of $i = 1 $ or $i = 2$.  In this case we 
will say that $(V_1, V_2)$ is the {\it induced Heegaard splitting of  $E(K)$} induced by $(V_1^i, V_2^i), i = 1,2$.  

\vskip20pt

\section{Interior tunnels}

\vskip10pt

Consider now a Heegaard splitting $(V_1 , V_2)$  for $E(K)$ the exterior of $K = K_1 \# K_2$, 
where  $\partial E(K) \subset V_1$ and  in which the decomposing annulus $A$ meets $V_1$ in disks 
and two vertical annuli.  Since the annulus $A$  meets $V_2$ in a connected planar surface $P$
it separates $V_2$ into two components each of which is a handlebody. We will denote the 
handlebodies $cl(V_2 - A) \cap E(K_i)$ by $V_2^i$ respectively. 
However $V_1 - A $ might have many components.

\begin{definition}
\label{floatdef}
\rm A component of $cl(V_1 - A) $ which is disjoint from $\partial E(K_i)$ and 
intersects $A$ in $n$ disks will be called an  {\it n-float} (see Fig. 2).
\end{definition}

\noindent {\bf Remark:}  Note that a n-float is either a 3-ball or a handlebody if its spine is 
not a tree. Furthermore there are always exactly two components of $cl(V_1 - A) $ not disjoint 
from $\partial E(K_i)$ (one in each of $E(K_1)$ and $E(K_2)$) and each one is a handlebody 
of genus at least one as $V_1$ is a compression body with a $T^2$ boundary. We denote these 
special  components by $N_1$ and $N_2$ depending on whether they are contained in 
$E(K_1)$ or $E(K_2)$ respectively.

\bigskip
Consider now any one of the meridian disks $E_i \subset \cal E$ of $V_2$. On $E_i$ we have a 
collection of arcs corresponding to the intersection with the decomposing annulus. These arcs, 
as indicated in Fig. 1, separate $E_i$ into sub-disks where disks on opposite sides of arcs  
are contained in opposite sides of $A$ i.e., in $E(K_1)$ or  $E(K_2)$ respectively. So each sub-disk 
is contained in either $E(K_1)$ or  $E(K_2)$. The boundary of these sub-disks is a collection 
of alternating arcs $\cup(\alpha_i  \cup \beta_i)$ where $\alpha_i$ are arcs on $A$ and $\beta_i$
are arcs on some component of  $cl(V_1 - A)$.

\vskip10pt

\vbox{\hskip70pt \epsfysize200pt\epsfbox{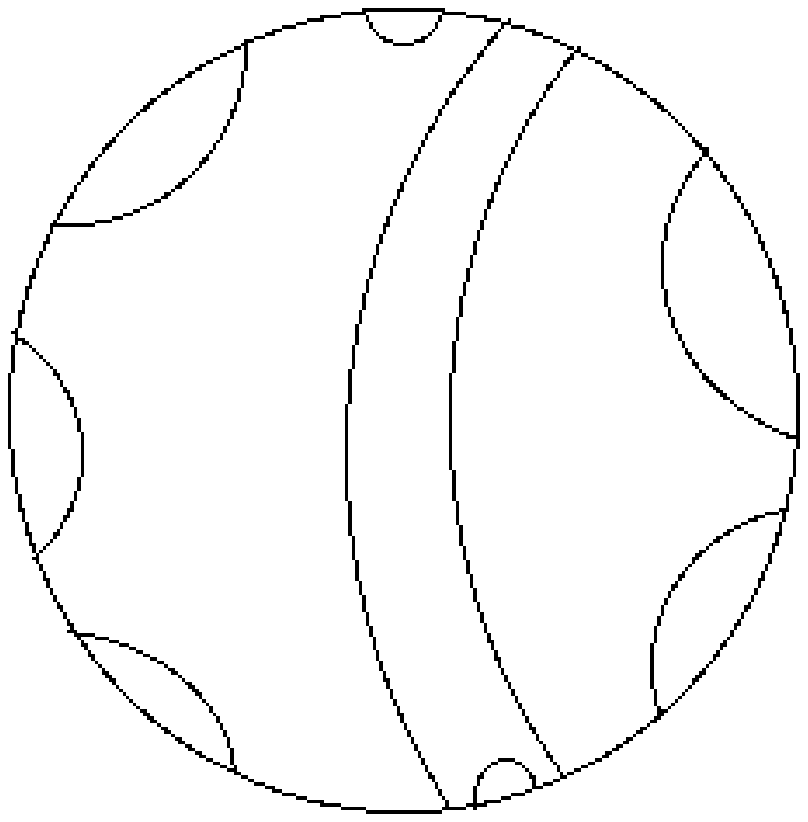}

\vskip10pt\centerline{Fig. 1}}

\vskip15pt

\begin{proposition}
\label{innertunnelProp}
Let $K_1$ and $K_2$ be knots in $S^3$ and let $K, A , \cal E$ be the  connected sum, a 
minimal intersection decomposing annulus and a meridional system for some Heegaard splitting
of $E(K)$ as above. Then
 
\begin{itemize} 
\item[\rm (a)] the $\beta$ arc part of the boundary of an outermost sub-disk in $E$ cannot
 be contained in a n-float of genus $0$.
\item [\rm (b)] if the $\beta$ arc part of the boundary of an outermost  sub-disk in $E$ is contained 
in an $N_i$ component, $i = 1$ or $2$,  and if $K_i, i = 1,2$ are prime the genus of $N_i$ is greater 
than  one.
\end{itemize}

\end{proposition}

\begin{proof}
Denote an outermost sub-disk of some $E_j$ by $ \Delta$  and suppose it is cut off by an arc $\alpha$ 
on $A$. By the \lq\lq Facts" proved in [Mo1] pp 41 - 42, any such outermost arc $\alpha $ must have both 
end points on a single disk $D_i$ which belongs to some n-float of genus $0$. Furthermore 
$\alpha \cup D_i \subset A$ must  separate the boundary components of $A$. Assume further that 
$\partial \Delta = \alpha \cup \beta$ where $\beta$ is an arc on the n-float meeting $D_i$ in exactly two 
points $\partial \beta = \partial \alpha$. On $\partial D_i$  there is a small arc $\gamma$ so that $\gamma 
\cup \beta$ is a simple closed curve  on the n-float bounding a disk $D$ there, since the n-float has no
genus (see Fig. 2 below). Furthermore $\gamma \cup \alpha$ is a simple closed curve on $A$ which together 
with a boundary component of $A$ bounds a sub-annulus  of $A$. Hence $\gamma \cup \alpha$  bounds a 
disk $D'$ on the decomposing $2$-sphere of $K$ intersecting $K$ in a single point. Thus we obtain a 
$2$-sphere $D \cup \Delta \cup D'$ which intersects the knot $K$ in a single point. This is a contradiction 
which finishes case (a). 

For case (b), assume that  the outermost disk $\Delta$ is contained in $N_1$, say, and that genus 
$N_1$ is one. As before we have $\partial \Delta = \alpha \cup \beta$ where $\beta$ is an arc on 
$N_1$ and a small arc $\gamma$ so that $ \gamma \cup \beta$ is a simple closed curve on $N_1$. 
If $ \gamma \cup \beta$ bounds a disk in $N_1$ we have the same proof as in case (a). If
 $ \gamma \cup \beta$ does not bound a disk on $N_1$ we consider small sub-arcs $\beta_1$
and $\beta_2$ of $\beta$ which are respective closed neighborhoods of $\partial \beta$. These
arcs together with a small arc $\delta$ on $\partial N_1  - \partial E(K_1)$  and $\gamma$ bound 
a small band $b$ on $\partial N_1$. Notice that $b  \cup_ {\beta_1,\beta_2}\Delta$ is an annulus
$A'$. The annulus $A'$ together with the sub-annulus $A''$ of $A$ cut off by $\alpha \cup \gamma$ 
defines an annulus $A' \cup_{\alpha \cup \gamma} A''$ which determines an isotopy of a meridian 
curve in $\partial A$ to a simple closed curve $\lambda$ on $\partial N_1$. Note that $N_1$ is a solid torus 
and $\pi_1(N_1) =  \integers$ which is generated by a meridian $\mu$ of $E(K_1)$. Hence 
$[\lambda]  = \mu \in \pi_1(N_1)$ (see Fig. 3).

Now we can consider the  
annulus $(A - A'') \cup A'$. If it is non-boundary parallel then since both knots $K_1, K_2$ are prime it must be
a  decomposing annulus which has at least one less disk component intersection than $A$ in contradiction to the
choice of $A$.  If it is boundary parallel, then as above, we have $A'' \cup A'$ as a decomposing annulus with a
smaller  number of disks. Again in contradiction  to the choice of $A$. So genus $N_1$ cannot be one and this
finishes case (b).

 \end{proof}

\vskip10pt

\begin {corollary}

\label{innertunnel}

Let $K_1, K_2 \subset S^3 $ be prime knots. Then every unknotting tunnel system $\tau$ for 
$K= K_1 \# K_2$ must contain at least one tunnel which is disjoint from the decomposing annulus 
for $K$ which minimizes the number of  intersections with  $N(K\cup \tau)$. 
\end{corollary}

\begin{proof}
Since the $\beta$ part of an outer-most disk must be contained in a float of genus greater than 
one we must have a $1-handle$ on the float to create the genus. The $1$-handle is the tunnel
disjoint from the decomposing annulus $A$.

\end{proof}

\vglue-60pt{\hskip50pt \epsfysize260pt\epsfbox{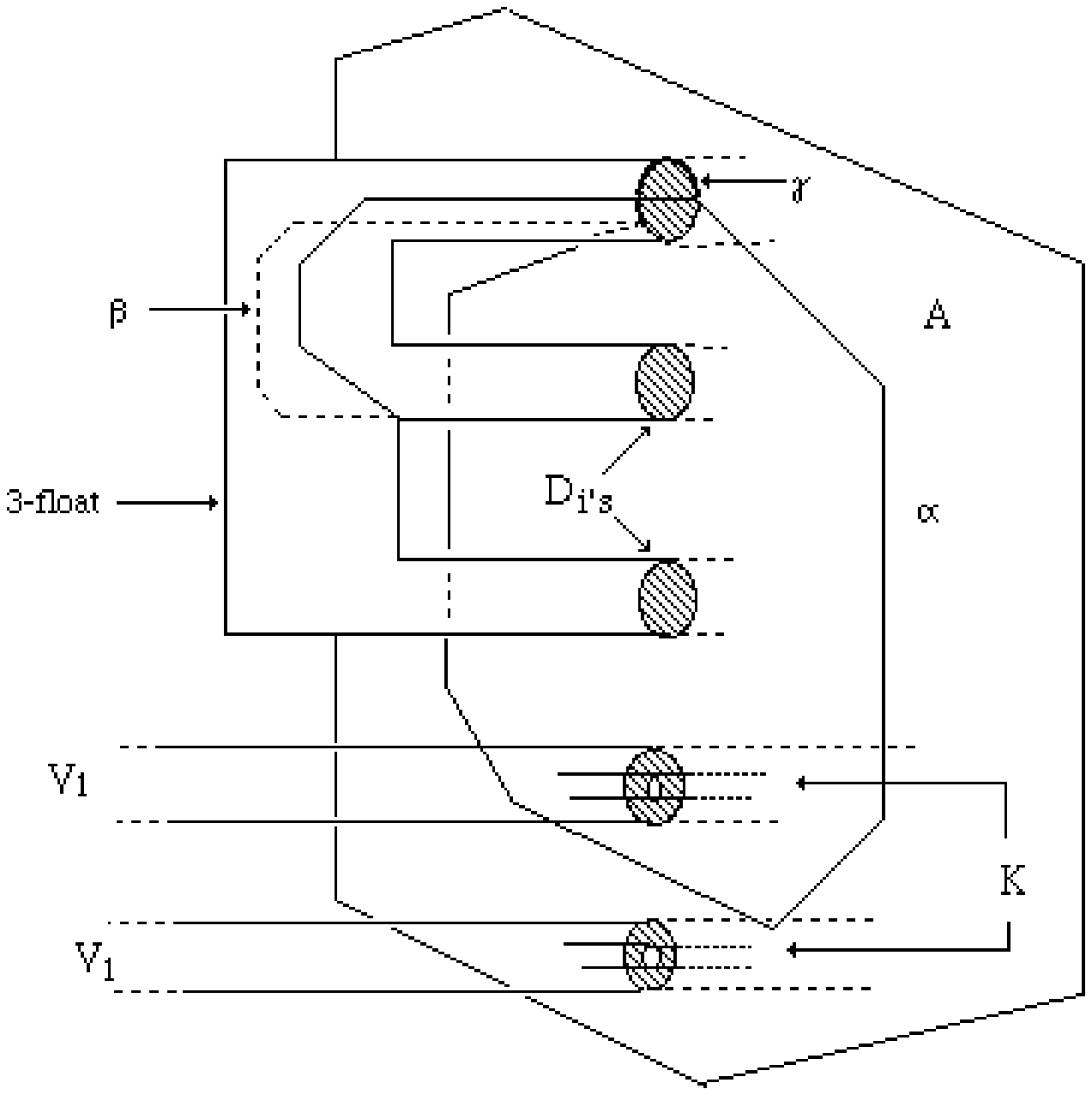}

\vskip15pt\centerline{Fig. 2}}

\vskip15pt

\vbox{\hskip50pt \epsfysize250pt\epsfbox{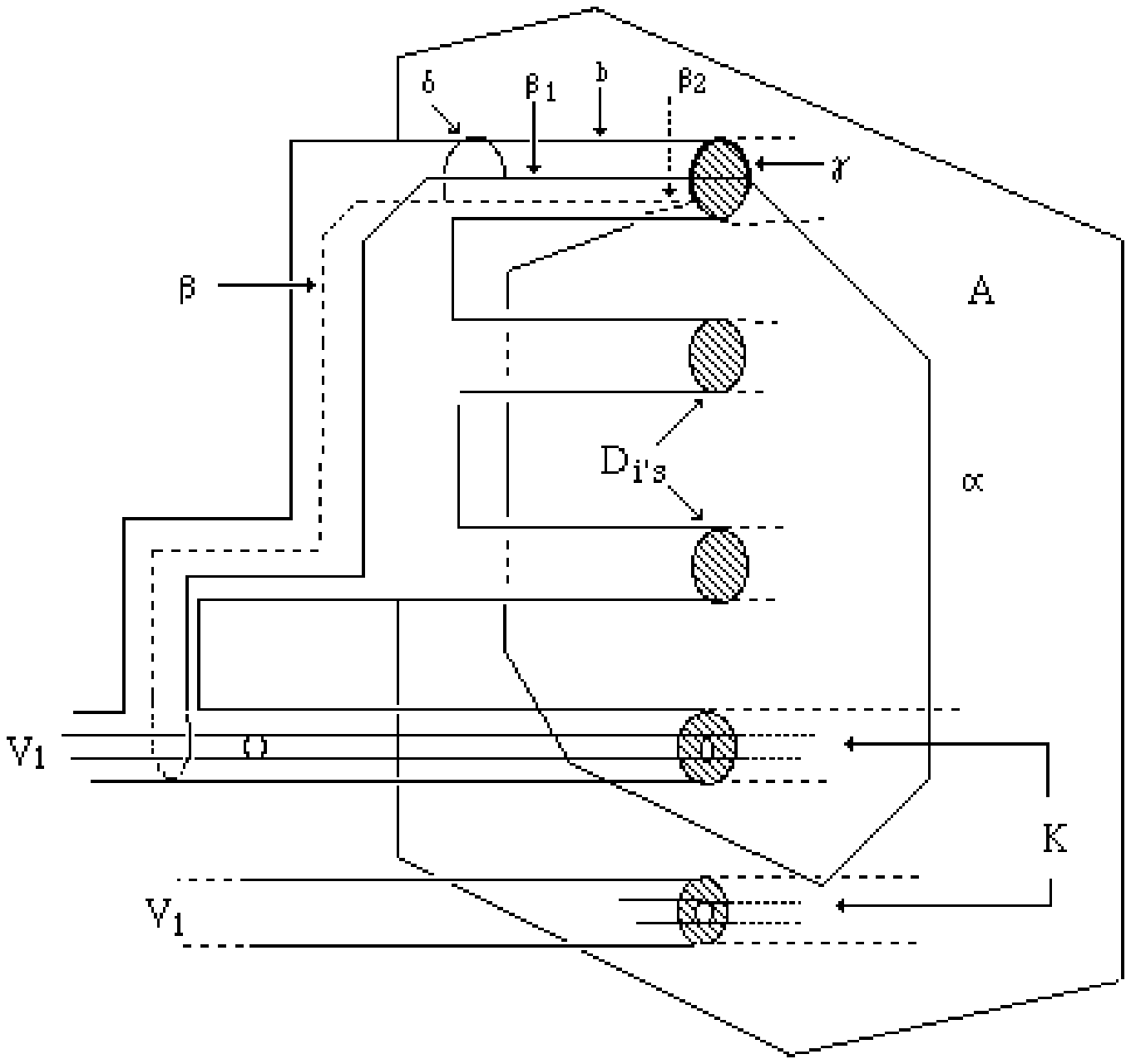}

\vskip10pt\centerline{Fig. 3}}

\vskip20pt

\section{Super additive and additive knots I}

%Weakly  irreducible Heegaard splittings I
\label{weak I}

\vskip 15pt

Given knots  $K_1, K_2 \subset S^3$ then the knot $K = K_1 \# K_2$ falls into one of three possibilities.

\begin{itemize}

\item[(i)] $t(K) = t(K_1) + t(K_2) + 1$ 

\item[(ii)]  $t(K) = t(K_1) + t(K_2) $

\item[(iii)] $t(K) \leq t(K_1) + t(K_2) - 1$

\end {itemize}

Recall that Heegaard splittings $(V_1^i , V_2^i ), i =1,2$ of $E(K_i)$ induce a Heegaard splitting 
$(V_1, V_2)$ of  $E(K_1 \# K_2)$ if and only if one of $(V_1^i , V_2^i ), i =1,2$  has a primitive 
meridian. If$(V_1, V_2)$ is induced then  we have $t(K) \leq t(K_1) + t(K_2) $. Therefore Case (ii) 
splits into two  subcases: (a) $(V_1, V_2)$  is induced by $(V_1^i , V_2^i), i =1,2$ and (b) $(V_1, V_2)$  is
not induced by $(V_1^i , V_2^i), i =1,2$. In this section we will deal with Case(i) and Case(ii) (a).
Case(ii) (b) will be discussed  in the next section.  In Case (i) we have:

\vskip13pt

\begin{theorem}
\label{Supaddthm}
Given knots $K_1, K_2$  and  $ K = K_1 \# K_2$ in $S^3$  for which the tunnel number 
satisfies $t(K) = t(K_1) + t(K_2) +1$ i.e., $t(K)$ is super additive,  then there is a minimal 
genus Heegaard splitting of  $E(K)$ which is  weakly reducible.
\end{theorem}

\vskip8pt

\begin{proof}
No one of the two knots has a Heegaard splitting where the meridian is a primitive
element since $t(K) = t(K_1) + t(K_2) +1$. A primitive meridian would mean that the Heegaard splittings 
of the knots will induce a Heegaard splitting of the connected sum which would make the tunnel number
 additive or less. Now drill a tunnel in $V_2^i$ with end points on opposite sides of the meridian curve on
$\partial V_2^i$ for one of the knots $K_i$ and add it as a $1$-tunnel to $V_1^i$ thus making the meridian 
primitive at the expense of increasing the genus by $1$. The two Heegaard splittings will now induce a 
Heegaard splitting on the connected sum which is of genus $t(K) + 1$. It is minimal since $t(K)  =  g - 1$ and 
weakly reducible by Proposition  \ref{wredProp}.

\end{proof}

\noindent For Case (ii) (a) we have the following theorem:

\begin{theorem}
\label {wredProp}
 Let $K_1, K_2$  and  $ K = K_1 \# K_2$ be knots  in $S^3$
and $(V_1^i , V_2^i ), i =1,2$ be  Heegaard splittings  for $E(K_i)$. If $(V_1^1, V_2^1)$  and $(V_1^2 , V_2^2)$
induce a  Heegaard splitting $(V_1, V_2)$  of  $E(K)$  then $(V_1, V_2)$ is a weakly reducible Heegaard splitting.
.

\end{theorem}

\begin{proof} 
We can assume that the decomposing annulus $A$ intersects 
the Heegaard splitting  $V_1, V_2$ as follows: It intersects $V_1$ in two vertical annuli 
and $V_2$ in one meridional annulus. (This is a consequence of the fact that $(V_1, V_2)$
is induced by the respective Heegaard splittings). Choose two essential disks $D_1^1$ and 
$D_1^2$ for $V_1$ on both sides of $A$, for example  cocore disks for  tunnels. Note that 
$D_1^1 \subset V_1^1$ and $D_1^2  \subset V_1^2$. The handlebody  $V_2$ is obtained 
from $V_2^1$  and $V_2^2$ by gluing them along the meridional annulus $A$. Since $V_2$  
turns out to be a handlebody $A$ must  be a primitive  annulus in at least one of $V_2^1$ or 
$V_2^2$, say $V_2^1$. So
there is at least one essential  disk $D_2$ in $V_2^1$ which is disjoint from $A$ and hence is 
also an essential  disk in $V_2$.  But $D_2$ is disjoint  from  $D_1^2$ as $D_1^2$ is also 
disjoint from $A$ and is on the opposite side.  Hence the Heegaard  splitting $(V_1, V_2)$ is 
weakly reducible.

\end{proof}

\vskip8pt

\begin{remark}\rm

Case (i) (a) is very common indeed e.g., any two knots which realize a minimal Heegaard splitting in a 
$2$-plat projection with the canonical tunnel systems will have a weakly reducible Heegaard splitting 
when composed (see [LM])

\end{remark}

\vskip25pt

\section{Additive knots II}

\vskip15pt

In this section we consider Case(ii) (b):  In this case both knots cannot have minimal genus Heegaard 
splittings with primitive meridians. Knots with this property, called also {\it fiendish knots}, are very elusive and 
their existence was first proved in [MR] and first examples were given in [MSY]. The knots considered 
in both [MR] and [MSY] satisfy $t(K) = t(K_1) + t(K_2) + 1$ so they fall into Case(i). For fiendish knots 
we have the following conjecture (see also Conjecture 1.5 of [Mo4]):

\begin{conjecture}
\label{sndconj}
Knots $K_1,K_2 \subset S^3$ will satisfy $t(K) = t(K_1)+t(K_2)+1$ 
if and only if both $E(K_1)$ and $E(K_2)$ do not have minimal genus Heegaard splittings with
primitive meridians.

\end{conjecture}

Note that Conjecture \ref{sndconj} implies Conjecture \ref{fstconj}. As if  Conjecture \ref{sndconj} is true 
then Case (ii)(b) cannot arise as all such knots will be in Case(i) and we are done. Conjecture \ref{sndconj} 
is known for knots which do not contain essential surfaces with meridian boundary components ([Mo4]
Theorem 1.6). 
We have the following:

\vskip12pt

\begin{definition}
\label{horDef}
\rm An incompressible meridional surface $S$ in a knot complement $E(K)$ will be called  {\it $\Sigma$  
horizontal} if it is not an annulus and it is contained in a Heegaard surface $\Sigma$ of $E(K)$ as a 
sub-surface, except for annuli collar neighborhoods of the meridian boundary components of $S$.  
These annuli will have one boundary component on the surface $\Sigma$ and the other on $\partial E(K)$.

\end{definition}

\vskip8pt

\begin{theorem}
\label{nhorwThm}
Let $K = K_1 \# K_2 \subset S^3$ be a knot. Any Heegaard surface $\Sigma$ for
$E(K)$ which does not contain any $\Sigma$ horizontal surfaces is weakly reducible.
\end{theorem}

\vskip8pt

\begin{proof}
Assume in contradiction that $(V_1, V_2)$  is a  strongly irreducible Heegaard 
splitting for  $E(K_1 \# K_2)$. Let $\Sigma = \partial V_1 = \partial V_2$ be the Heegaard surface
and let $A$ be the decomposing annulus for the connected sum minimizing the intersection with 
$\Sigma$. We can assume (see Lemma 2.3 of [Mo3]) that after an isotopy of the annulus $A \cap \Sigma$
is a collection of essential curves on both $A$ and $\Sigma$. Hence, as we assumed that $V_1$ is
the compression body containing $\partial E(K_1 \# K_2)$ then $V_1 \cap A$ is composed of
two vertical annuli $A_1^*,  A_2^*$ and a minimal collection of essential annuli $A_1, \dots,  A_d$ and
$V_2 \cap A$ is composed of a minimal collection of essential annuli $B_1, \dots, B_{d + 1}$. By Lemma 
2.1 of  [Mo3] we can find essential disks $D_1, D_2$  in $V_1, V_2$ respectively which are disjoint
from $A_1, \dots,  A_d$ and $B_1, \dots, B_{d + 1}$.  Since  $A_1^*,  A_2^*$ share a boundary
component with $B_1$  and $B_{d + 1}$ we can conclude that the disks  $D_1, D_2$ are disjoint
from $A$. The annulus $A$ splits each of $V_1$ and $V_2$ into two unions of handlebodies $\cup_r V_{1,r}^i$
and $\cup_s V_{2,s}^i$  respectively where i = 1,2 depending if the component is in $E(K_1)$ or $E(K_2)$. 
If the disks $D_1, D_2$ are contained in $V_{1,r}^i$ and $  V_{2,s}^j$ respectively, $i, j \in \{ 1, 2\}$,
for different values of $i$ and $j$ then  $\partial D_1 \cap \partial D_2 = \emptyset$ as both
of $D_1$ and  $D_2$ are disjoint from $A$. Hence the Heegaard splitting $(V_1, V_2)$ is weakly
reducible in contradiction. So we can assume that both of $D_1$ and  $D_2$ are contained in
$V_{1,r}^i$ and $V_{2,s}^i$ for the same $i$, say $i =1$ i.e., on the same side of $A$. Consider  now the
components of
$\Sigma -A$ contained in  $V_1^2$ and $ V_2^2$. An innermost disk argument shows that each of these 
components must be incompressible in $V_1^2$ and $ V_2^2$ as otherwise we obtain a 
compressing disk $D_3$ disjoint from $A$  which is disjoint from both $D_1$ and  $D_2$ 
and hence the Heegaard splitting  $(V_1, V_2)$  is weakly reducible in contradiction. The boundary
curves of any component of  $\Sigma -A$ contained in  $V_1^2$ and $ V_2^2$ are essential curves 
on the meridional decomposing annulus $A$ and hence are isotopic to meridian curves in $E(K_2)$.
Therefore they are isotopic to meridian curves in $E(K)$. Thus these components of $\Sigma -A$ are
horizontal surfaces. Since we assumed that such surfaces do not exist in  $E(K)$  we obtain a 
contradiction to our  assumption that $(V_1, V_2)$ is a strongly irreducible Heegaard splitting 
of $E(K)$.

\end{proof}

\vskip10pt

\begin{remark}
\label{Mremark}
\rm   A result of similar nature is mentioned by  Morimoto (see [Mo3]  Remark 4.3): If 
$K_i \subset M_i$ are knots  then $E(K_1 \# K_2)$  always has a weakly  reducible 
Heegaard splitting of minimal genus if  none of  $M_1$ and $M_2$ have Lens space summands and 
none of $E(K_1)$ and $E(K_2)$ contains meridional essential surfaces. It seems that the conditions in
Theorem  \ref {nhorwThm} are weaker. 
\end{remark}

\vskip10pt

We will now specialized to the situation where there is a tunnel system for $K$ with a single tunnel 
minimally intersecting the decomposing annulus in a single point. More precisely: $E(K)$ has a minimal genus 
Heegaard splitting so that $t(K) = t(K_1) + t(K_2)$ and $V_1 \cap A$ consists of two spanning annuli 
and a single disk. This is clearly a subset of Case (ii) (b). However to the best of my knowledge all examples 
of minimal tunnels systems of composite knots which have tunnels intersecting the decomposing annulus 
essentially do so exactly once. 

Before we specialize we need the theorem below which is true in a more general setting. It  is of independent 
interest as it gives a new  characterization for when a set of primitive curves on a handlebody is simultaneously 
primitive (compare [Go]). 

Given a collection of annuli $A_1, \dots, A_n$ on the boundary of a handlebody 
$H$ we say that they are {\it simultaneously primitive} if there exists a  collection  $D_1, \dots, D_n$ 
of disjoint essential disks so that $D_i \cap A_i$ is a single essential arc in $A_i$  and if $i \neq j$ then 
$D_i \cap A_j = \emptyset$. 

\vskip 30pt

\begin{theorem}
\label{simprim}
Let $H_1$ and $H_2$ be two handlebodies and let  $B_1,\dots, B_n$ be a set of disjoint non-parallel 
incompressible primitive annuli in $\partial H_1$.  Let $C_1,\dots, C_n$ be any collection  of 
incompressible \underline  {non} primitive  disjoint annuli in $\partial H_2$. Then  $B_1,\dots, B_n$ are 
simultaneously primitive in $H_1$ if and only if  $H_1 \cup_{\{B_1=C_1, \dots, B_n=C_n\}}  H_2$ 
is a handlebody. 

\end{theorem}

\vskip10pt

\begin{proof}

Assume first that the annuli $B_1,\dots, B_n$ are simultaneously primitive in $H_1$.
The proof will be by induction on $n$.  For $n = 1$ we can glue $H_1$ to $H_2$ along 
$B_1$ and $C_1$ to obtain a manifold $N_1$. Since  the annuli $B_1$ and $C_1$ are 
incompressible we have that $\pi_1(N_1) = \pi_1(H_1)  *_{\integers} \pi_1(H_2)$. The 
generator of the $\integers$ is a primitive element in the free  group $ \pi_1(H_1)$ so $\pi_1(N_1)$ 
is a free group. It now follows from the Loop Theorem that $N_1$ is a handlebody.  Assume 
by induction that $N_{n-1} = H_1 \cup_{\{B_1=C_1, \dots, B_{n-1} = C_{n-1}\}}  H_2$ 
is a handlebody. The annulus $B_n$ is disjoint from the annuli $B_1, \dots, B_{n-1}$ and  
$C_1, \dots, C_n$ and is still primitive in $N_{n-1}$ as the annuli $B_1, \dots, B_n$ are 
simultaneously primitive and non-parallel and hence there is an essential  disk $D$ in $N_{n-1}$ 
which is disjoint from $B_1, \dots, B_{n-1}$ and  $C_1, \dots, C_n$ and which intersects $B_n$ 
in a single arc. Now $N_n$ is obtained from  $N_{n-1}$ by gluing the primitive annulus $B_n$ 
to the annulus $C_n$. Hence $\pi_1(N_n) =  \pi_1(H_{n-1})  *_{\integers}$ is an HNN extension 
of the free group $\pi_1(H_{n-1})$ where two $\integers$-subgroups are identified and the generator 
of one of them is a primitive element. It follows that $\pi_1(N_n)$ is a free group and again by the 
Loop Theorem $N_n = H_1  \cup_{\{B_1 = C_1, \dots, B_n = C_n \}}  H_2$ is a handlebody.

For the proof in the other direction: Assume that $H_1  \cup_{\{B_1 = C_1, \dots, B_n = C_n \}}  H_2$ 
is a handlebody $H$ and let $B = \{B_1, ..., B_n\}$ and $C = \{C_1, ..., C_n\}$ be as in the
theorem.  Let $H'_1$ be the result of cutting $H_1$ along a maximal set of compression 
disks of $\partial H_1 - \cup B_i$.  Note that gluing $H'_1$ to $H_2$ along $B$ and $C$ 
yields a handlebody. As it is obtained from the handlebody $H$ by cutting it along disks 
which are disjoint from both of $B$ and $C$. Up to relabeling we may assume that 
$B' = \{B_1, ..., B_k\}$ is the set of annuli in $B$ which are a longitudinal annulus of some solid 
torus component $V_i$, $i = 1, \dots, k$ of $H'_1$ containing no other $B_j$.  Denote by 
$H''_1 = H'_1 - \cup V_i$, and let $B'' = B - B'$.  There is no compressing disk in $H''_1$ 
intersecting $B''$ in a single essential arc. As any such disk would define another torus components
$V_j$ containing some annulus $B_j$, $j \notin \{1, \dots, k \}$ and no other annulus.

Let $C'$ and $C''$ be the corresponding subsets of $C$.  Then $H'_1
\cup _{B=C} H_2$ can be obtained by gluing $V_1, ..., V_k$ to $H_2$
along $B'$ and $C'$ to obtain a manifold $H'_2$, and then gluing
$H''_1$ to $H'_2$ along $B''$ and $C''$.  The manifold $H'_2$ is a handlebody and
is  homeomorphic to $H_2$ by the definition of $B'$ and the first part of the theorem.
Hence the annuli $C''$  are still non-primitive annuli on $\partial H'_2$.  If $B$ is not 
simultaneously primitive then $B''$ is non-empty, hence after gluing the remaining
components of $H'_1$ to $H'_2$, the surface $B'' = C''$ is an
essential surface in the handlebody $H'_1 \cup H_2 = H''_1 \cup H'_2$
because there is no compressing or boundary compressing disk for this
surface, which contradicts the fact that there are no essential
non-disk surfaces in a handlebody.   

\end{proof}
\vskip10pt

\noindent Further evidence in the direction of Conjecture \ref{fstconj} is the following:

\vskip 10pt

\begin{theorem}
\label{singtunthm}
 Let $K_1, K_2$ be prime knots in $S^3$ and  $K = K_1 \# K_2$. Assume that $t(K) = t(K_1) + t(K_2)$
and  $t(K_i) \leq 2 $.  Furthermore, assume that a minimal tunnel system for $K$ minimaly intersects a  
decomposing annulus $A$ in a single point, then there is a Heegaard splitting of  $E(K)$  of minimal genus 
which is weakly reducible.

\end{theorem}

\begin{proof} Let $(V_1, V_2)$ be the Heegaard splitting of $E(K)$ determined by the minimal tunnel 
system which intersects the decomposing annulus $A$ in a single point.  We can therefore assume that 
$V_1 \cap A = A_1^*\cup A_2^* \cup D_1$.  The once punctured annulus $A \cap V_2$ has two  boundary 
components coming from the vertical annuli  $A^*_1, A^*_2$ and denoted by $C^*_1, C^*_2$ 
respectively and one boundary component $\partial D_1$ coming from the tunnel. As $A$ intersects $V_1$ 
minimally $A - D_1$ is an incompressible planar surface in a handlebody and hence is boundary compressible.
A boundary compression cannot be on an arc connecting $C^*_i, i = 1,2$ to $\partial D_1$ as then we could
use the compressing disk to isotope the tunnel off $A$. Such an arc will be called of {\em Type I}.  Furthermore
a boundary compression cannot be on an arc connecting $C^*_1$ to $C^*_2$ as then $A$ will be boundary
parallel in contradiction. Hence the boundary compressing arc will connect $\partial D_1$  to itself and since it is 
non-trivial it must separate $C^*_1$ and $C^*_2$. Such an arc will be called of {\em Type II}.

Choose a meridional system of disks  ${ \cal E}  = E_1, \dots, E_{t(K) + 1}$ for $V_2$. Each disk
in $\cal E$ must intersect $D_1$ as otherwise the Heegaard splitting will be weakly reducible and we  are 
done. An outermost arc of intersection $\alpha$ on some $E_i$  separates a boundary compressing sub-disk 
$\Delta \subset E_i$ and from the previous paragraph  $\alpha$ is an arc of type II on $A$

\vskip 10pt

We can boundary compress $A$ along $\Delta$ or alternatively isotope $\partial V_1 = \partial V_2$ 
along $\Delta$. Doing the second operation  does not change $A$ or the isotopy class of the Heegaard 
splitting $(V_1, V_2)$, but does change the intersection of the \lq\lq new" Heegaard surface, also denoted by 
$\partial V_1 = \partial V_2$, with $A$. The result is that now  $A \cap V_1 = A_1^* \cup A_2^* \cup A_1$, 
where $A_1$ is an essential sub-annulus of $A$ which contains the disk $D_1$. The intersection 
$A\cap V_2 = B_1 \cup B_2$, where $B_1, B_2$ are also essential sub-annuli of $A$ (as in  Fig. 4).

%\vskip -0pt

\vbox{\hskip60pt \epsfysize210pt\epsfbox{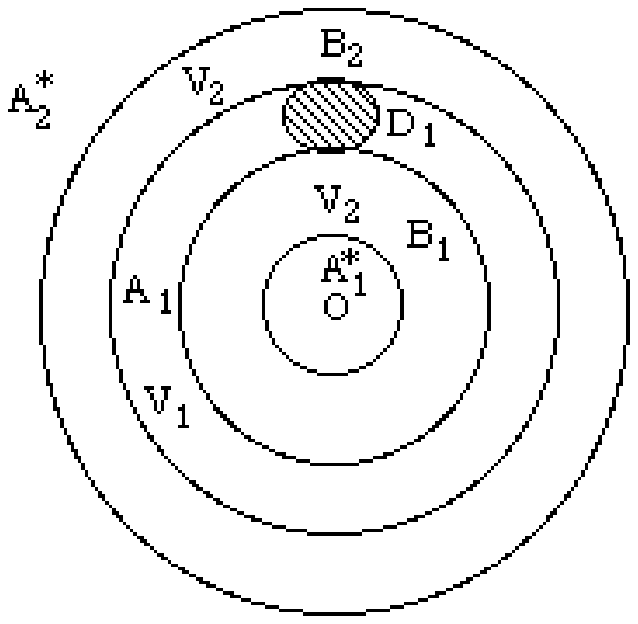}

\vskip10pt\centerline{Fig. 4}}

\vskip15pt

Let $V_1^i$  denote the components of $V_1 - A$ and $V_2^i$ denote the components of $V_2 - A$.  Assume that 
the disk $\Delta$ is contained in $E(K_k), k = 1$ or $k = 2$. Note that  isotoping the Heegaard surface $\partial V_1$ 
along $\Delta$ changes the \underline{induced} Heegaard splitting only on the knot complement containing 
$\Delta$, (i.e., on $E(K_k)$ only!).  On the induced Heegaard splitting of $E(K_k)$ this isotopy is equivalent 
to cutting ${V}_2^k$ along $\Delta$ to obtain $W_2^k$ and  adding the $2$-handle  $N(\Delta)$ to  ${V}_1^k$ 
to obtain $W_1^k$. It is possible that in this case $W_1^k$ might not be a handlebody. It is also possible that $\Delta$ is
a separating disk in  $V_2^k$ and in this case, $W_2^k$  might have two components  $W_2^{k,1}$ and $W_2^{k,2}$.

The annuli $B_1, B_2$ are essential annuli contained in $V_2$ which together separate $V_2$.
Hence, when we cut $V_2$ along them we obtain a handlebody $W_2^j, j \neq k$, and if neither 
of $B_1$ or $ B_2$ is separating a handlebody $W_2^k$. If one of  $B_1$ or $B_2$ is separating then 
$V_2 \cap E(K_k)$ splits into two handlebodies $W_2^{k,1}$ and $W_2^{k,2}$ 
This is  the situation corresponding to the disk $\Delta$ being a separating disk in ${V}_2^k$.
Denote the \lq\lq traces"  of $B_1$ and $B_2$ on $W_2^i$ by $B_1^i, B_2^i, \hskip3pt i = 1, 2$.

Since $t(K) = t(K_1) + t(K_2)$  and only one 
tunnel gets split into two arcs by cutting along $A$ it follows that after cutting $(V_1, V_2)$ 
along $A$ there are two possibilities: The induced  Heegaard splitting $(V_1^1, V_2^1)$ of $E(K_1)$
is of minimal genus  and $(V_1^2, V_2^2)$ of $E(K_2)$ is of minimal genus plus one or vice versa
Up to relabeling the knots we assume that $E(K_1)$ is of minimal genus.
\vskip10pt

\noindent Claim (1): If one of $B_1^1$ or $B_2^1$ is primitive in $W_2^1 \subset E(K_1)$ then 
$E(K)$ has a weakly reducible Heegaard splitting of minimal genus.

\vskip 10pt

\noindent Proof (Claim (1)): If the disk $\Delta $ is contained in $E(K_2)$ then $(W_1^1, W_2^1)$ is a
Heegaard splitting of minimal genus for $E(K_1)$. So, if either  $B_1^1$ or $B_2^1$ is a primitive annulus on 
$W_2^1$ (which ,in this case, is equal to $V_2^1$ less a collar ) we will treat an isotopic image of the primitive annulus
$B_1^1$ or $B_2^1$  respectively on $\partial E(K_1)$  as a decomposing annulus. Now glue $E(K_1)$ to
$E(K_2)$  along this annulus to obtain  a Heegaard splitting of $E(K)$ which is of  minimal  genus 
(as that of  $(V_1, V_2)$) and is weakly reducible by Theorem \ref{wredProp}.

If, on the other hand, the disk $\Delta$ is contained in $E(K_1)$ then recall that we obtain $V_2^1$ from 
$W_2^1$ by identifying together the two \lq\lq traces" (copies) of the disk $\Delta$ on $W_2^1$, 
i.e., adding a $1$-handle to these traces. These traces intersect both of $\partial B_1^1$ and  $\partial B_2^1$ in 
a single arc each. Hence if  one of $B_1^1$ or $B_2^1$ is primitive in $W_2^1$ it would also be primitive in
$V_2^1$,  regardless of whether $B_1^1$ and  $B_2^1$  are separating or not. We now use the same argument as
above  to obtain a weakly reducible Heegaard splitting of the same genus  as that of $(V_1,V_2)$ of $E(K)$.

\qed

\vskip8pt

Thus we can assume that both of the annuli $B_1^1$ and $B_2^1$ are not primitive in 
$W_2^1 \subset E(K_1)$.  Since $V_2 =W_2^2 \cup_{B_1^2 =  B_1^1, B_2^2  = B_2^1} W_2^1$ 
is a handlebody it follows that $B_1^2$ and $B_2^2$ must be primitive in $ W_2^2 $: 
By setting $B_1 = B_1^2$, $B_2 = B_2^2$ and $C_1 = B_1^1$, $C_2 = B_2^1$ we satisfy the 
conditions of  Theorem \ref{simprim} and can conclude that $B_1^2$  and $B_2^2$ are simultaneously 
primitive in $W_2^2$. If  it happens that $W_2^2$ has more than one component we certainly have disjoint annuli
intersecting disjoint disks in a single arc. We will refer to this situation as the annuli being {\it extended 
simultaneously primitive}.

\vskip8pt

\noindent Claim (2): If $B_1^2, B_2^2$ are simultaneously primitive or extended simultaneously  
primitive  on $W_2^2 \subset E(K_2)$ the complement with the non-minimal genus Heegaard splitting, 
then $E(K)$ is a weakly reducible Heegaard splitting of minimal genus.

\vskip 20pt
\noindent Proof (Claim (2)): The induced Heegaard splitting of $E(K_2)$ is of genus at least three since 
it is induced by a tunnel  system containing at least two tunnels, i.e., one interior tunnel (by Corollary 
\ref{innertunnel}) and the  \lq\lq half" tunnel coming from  the split tunnel crossing $A$.  
Assume that the disk $\Delta$ is contained in $E(K_2)$ so after cutting $V_2^2$ along 
$\Delta$ we obtain either a handlebody of genus at least two with two simultaneously primitive 
annuli on it or a disjoint union of two handebodies one of which has at least genus two with two 
extended simultaneously primitive annuli  on them.

Thus in both cases there is at least one essential disk $D_2$ in  $W_2^2$, (a separating disk in the first
case),  which is disjoint from  $B_1^2$ and $B_2^2$ and hence from $A$. 
Since $V_2 = W_2^2 \cup_{B_1^2 =  B_1^1, B_2^2  = B_2^1} W_2^1$, as before, the disk 
$D_2$  is an essential disk in $V_2$ which is disjoint from $A$ and hence from the essential disk 
$D_1^* \subset V_1 $ which is the image of the disk $D_1$ pushed slightly into $E(K_1)$. Thus the 
Heegaard splitting $(V_1, V_2)$ of $E(K)$ is weakly reducible and we are done (see Fig. 5).

%\vskip-75pt

\vbox{\hskip40pt \epsfysize200pt\epsfbox{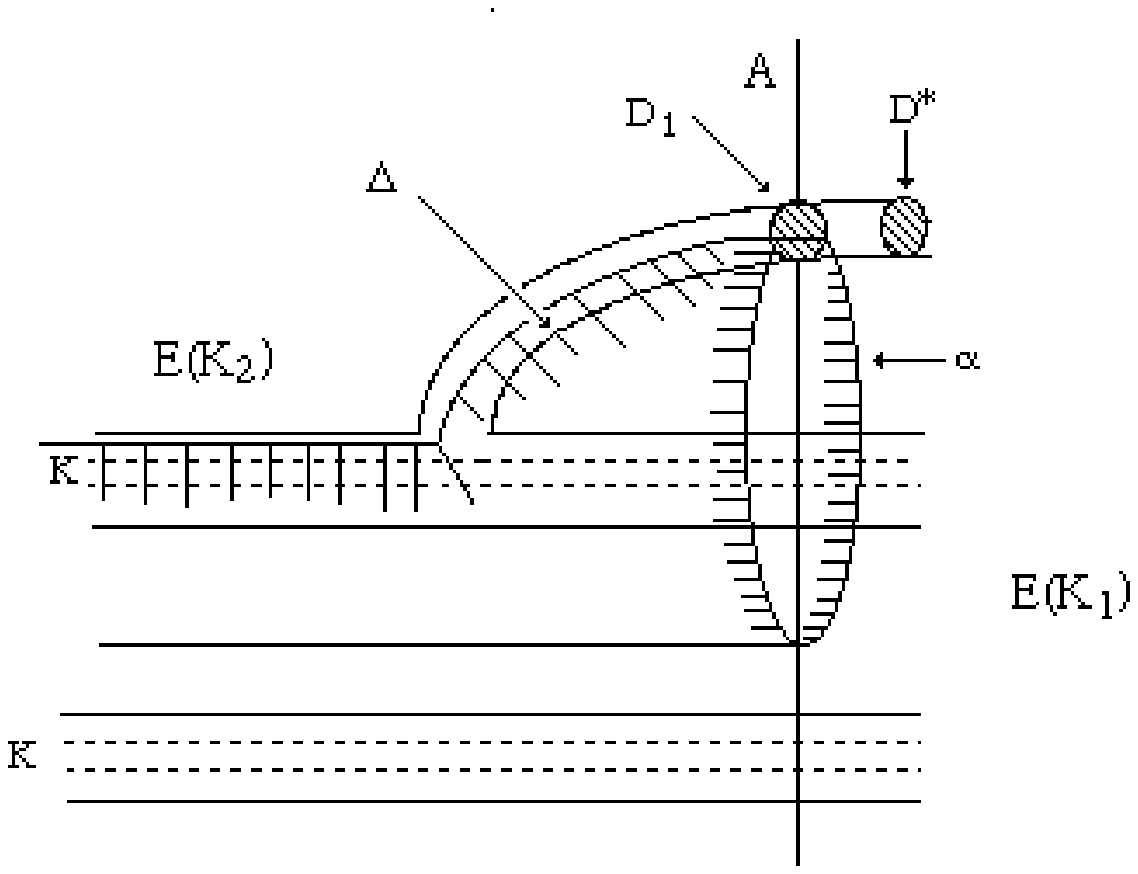}

\vskip10pt\centerline{Fig. 5}}

\vskip10pt

Assume therefore that the disk $\Delta$ is contained in $E(K_1)$. If $t(K_1) = 1$ then $V_2^2$
is a genus two handlebody and after cutting $V_2^1$ along $\Delta$ we obtain either one or two 
solid tori (depending if $\Delta$ is separating or not) embedded in $S^3$ with non-primitive annuli 
on their boundary. Extend these annuli all the way to $\partial E(K_1)$. When attaching disks to these 
meridional annuli one obtains a Lens space contained  in $S^3$, which is a contradiction.

If $t(K_1) = 2$ then $V_2^1$ is a genus three handlebody and after cutting $V_2^2$ along $\Delta$ we 
obtain either one solid torus component with one or two non-primitive annuli on its boundary (if $\Delta$
is separating) or a genus two handlebody with two non-primitive annuli on its boundary (if $\Delta$ is 
non-separating). In both cases the components are embedded in $S^3$. The first case is dealt with as in the 
previous paragraph. In the second case, then after adding two meridional disks along these annuli we obtain
a $2$-sphere $S \subset (S^3,K_1)$ which intersects $K_1$ in four points. In particular $S$ bounds a $3$-ball on 
both sides. If we change the order of cutting along $\Delta$ and adding disks by first adding the two meridional disks 
to the meridional annuli on  $V_2^1$ we obtain a solid torus $W_2$ with $\Delta$ as its unique meridional disk.
Since adding $N(\Delta)$ to the complement of $W_2$ gives us a $3$-ball then $S^3 - W_2$ is also a solid torus
$W_1$. The solid  torus $W_1$ can also be obtained from the genus three  compression body $V_1^1$ as follows: 
Fill  $\partial_{-}V_1^1$ with $N(K_1)$ to get a pair $( V, K_1)$. Now cut the pair $( V, K_1)$ along 
meridional disk corresponding to the meridional annul.i on $\partial V_1^1$. These annuli are not parallel 
on $\partial V_1^1$ so we get a solid torus $W_1$ whose unique meridian disk $\Delta'$ is a cocore disk of one of
the $1$-handles of $V_1^1$.  

Now since we obtained $S$ from $W_1$ and $W_2$ the disks $\Delta$  and $\Delta'$ are a canceling pair. 
But this implies that  the minimal genus Heegaard splitting $(V_1^1, V_2^1)$ is reducible in contradiction.
Hence this case cannot happen and the proof is complete.

\end{proof}

\vskip30pt

\section{Sub-additive knots}

\vskip10pt

In this section we consider connected sums of knots $K_n  \subset S^3$ as in Fig. 6 and  $2$-bridge knots 
$K({\frac \alpha \beta}) \subset S^3$  determined by $ {\frac \alpha  \beta} \subset \rationals $. 
These are the only examples so far of prime knots $K_1, K_2 \subset S^3$ so that $t(K) = t(K_1) + t(K_2) -1$. 
For these examples we have:

\begin{theorem}
\label{strthm} 
Let  $K_n \subset S^3$ be the knot as in Fig.6 and  $K({\frac \alpha \beta})  \subset S^3$ a 2-bridge knot
determined by $ {\frac \alpha \beta} \subset \rationals $. Let $K$ denote the connected sum 
$K_n \# K({\frac \alpha \beta})$, then the Heegaard splitting of $E(K)$ determined by the minimal tunnel
system  for $K$, (as in Fig. 6) is strongly irreducible.
\end{theorem}

%\vglue-170pt

\vbox{\hskip30pt \epsfysize150pt\epsfbox{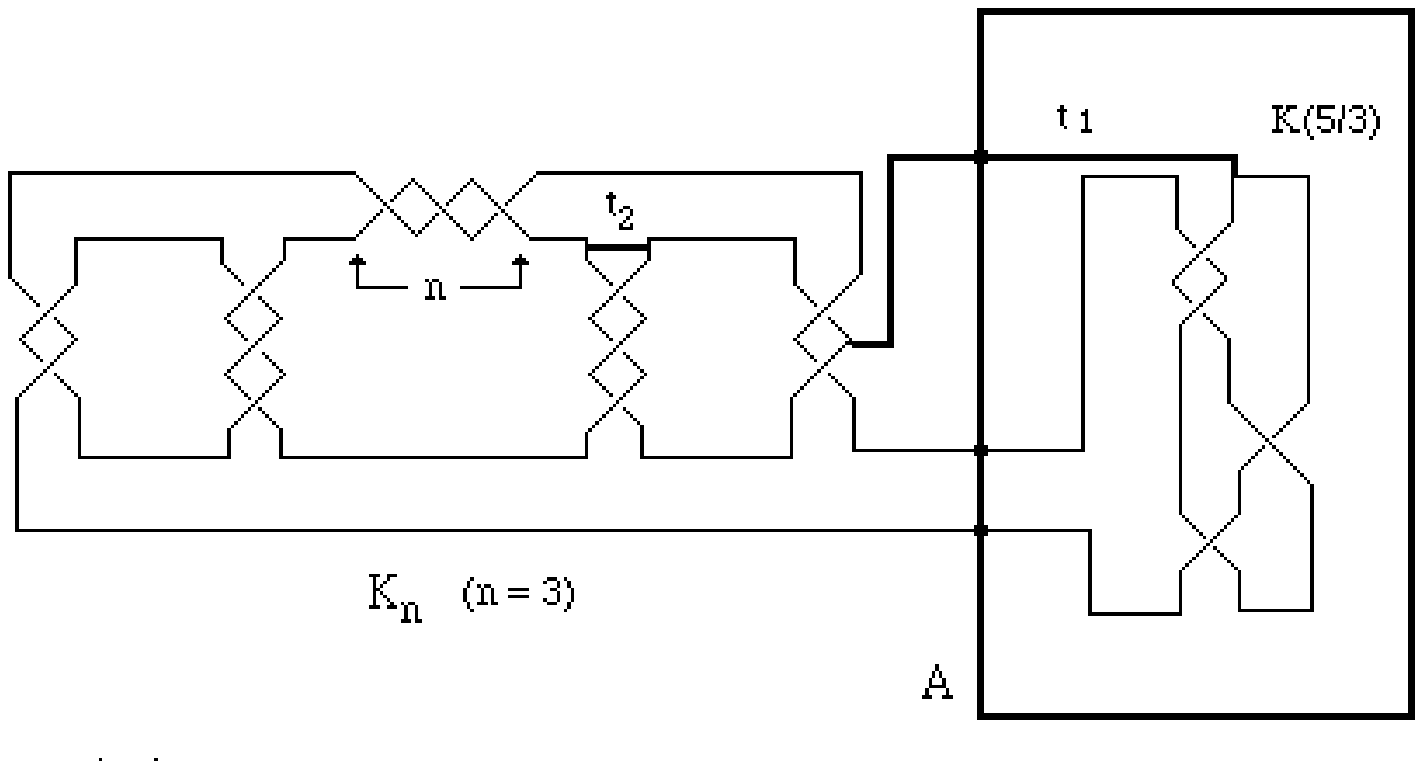}

\vskip8pt\centerline{Fig. 6}}

\vskip15pt

In  Fig. 6, $A$ denotes the decomposing annulus and $t_1$,$t_2$ denote the unknotting tunnels.

\vskip 10pt

\begin{proof}  Since $E( K({\frac \alpha \beta}))$  has a genus two Heegaard splitting
(as $ K({\frac \alpha \beta})$ is a tunnel number one knot) and is irreducible, the Heegaard splitting
is strongly irreducible. Otherwise we could compress the Heegaard surface to both sides and obtain an
essential $2$-sphere in contradiction. Similarly any Heegaard splitting of minimal genus three of a hyperbolic 
knot is strongly irreducible: As the knot complement is irreducible we can compress at most twice (once to
each side). But then, by compressing the Heegaard surface we obtain an incompressible non-boundary
parallel torus in contradiction to the fact that the knot is hyperbolic (see [Mh1]). The knots $K_n$ 
are alternating knots and not torus knots so by  Corollary 1 of [Me] they do no contain  incompressible 
non-boundary parallel tori and hence are hyperbolic. 

Note that $E(K)$ induces minimal genus Heegaard splittings, of genus two and three respectively, 
on both of $E( K({\frac \alpha\beta}))$ and $E(K_n)$. By slightly abusing notation we will denote 
the components of $E(K) - A$ by $E( K({\frac \alpha \beta}))$ and $E(K_n)$.

As in Fig.7 let $D$ denote the cocore disk of the tunnel
$t_1$ which   intersects the decomposing annulus $A$ and let $D'$ denote the cocore disk of $t_2$ the tunnel
interior  to $E(K_n)$. We can choose the disks ${\cal E} = \{D, D' \}$ as a meridional system of disks for the 
compression body  $V_1$. Note also that  $A$ minimizes the intersection with $V_1$ as if 
$A \cap V_1 = \emptyset$ the Heegaard genus of $E(K)$ would be additive and equal to three.
 
Let $F$ be the Heegaard splitting surface $\partial V_1 = \partial V_2$, 
and let $F_1 = F \cap E(K(\frac{\alpha}{\beta}))$,  and $F_2 = F \cap E(K_n)$.  For each
essential disk ${\cal D}_1$, ${\cal D}_2$ in $V_1$ and $V_2$ respectively, we choose a representative
in their isotopy class so that ${\cal D}_i \cap A$ is minimal; in particular, each component of 
$\partial {\cal D}_i \cap F_j$, $i, j \in \{1,2\}$, is an essential circle or essential arc on $F_j$, and each
component of ${\cal D}_i \cap A$ is an arc.

\vskip 15pt

\noindent Claim : Let $ E$ be an essential disk in $V_1$  then:

\begin{itemize}

\item [(a)] If $E \cap A \neq \emptyset$ then the outermost sub-disk  $E^\#$ of $E - A$  
is an essential disk in the components $V_1^2 \subset  E(K_n)$ or  $V_1^1 \subset E( K({\frac \alpha
\beta}))$ of $V_1 - A $,  depending on which side of $A$ contains $E^\#$. If it is in 
$E( K({\frac \alpha \beta}))$ then  $\partial E^\# = \gamma \cup \delta$ where $\gamma$ 
is an inessential arc on one of the vertical annuli  $A^*_i$ and  $\delta$ is an arc on 
$\partial V_1 - A$ as indicated in Fig. 7 .

\item [(b)] If $E \cap A = \emptyset$ and $E$ is contained in the $E( K({\frac \alpha \beta}))$ 
component then  $E$ is parallel to $D$.

\end{itemize}

%\vglue-70pt

\vbox{\hskip30pt \epsfysize230pt\epsfbox{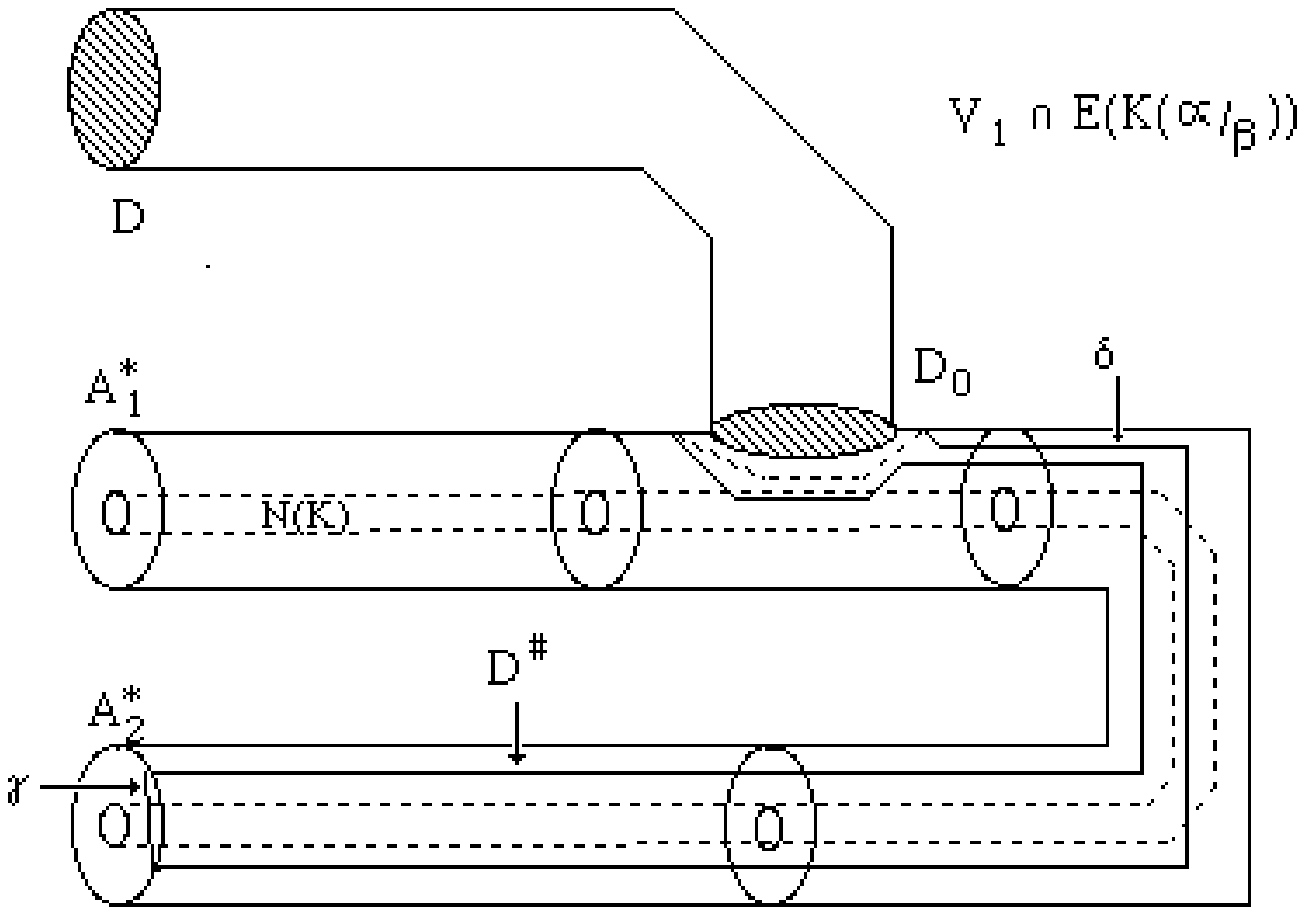}

\vskip15pt\centerline{Fig. 7}}

\vskip15pt

\begin{proof}  a)  Note that $\partial E^\#$ is the union of two arcs
$\gamma \subset A$ and $\delta$. If $E^\#$ is inessential we could isotope $E^\#$ 
off  $A$.  This is a contradiction to the choice of $E$.

Assume now that $E^\#$  is contained in  $E( K({\frac \alpha \beta}))$. If  $\gamma \subset D$
then $\partial E^\#$ is isotopic to a curve which represents a power of the meridian in $\pi_1(E(K))$ 
which is a contradiction as the meridian has infinite order in $\pi_1(E(K))$. So  $E^\#  \cap D = \emptyset$.
Consider now the disk $D_0$ which is the  intersection of $N(t_1)$ with  the component of $N(\partial E(K)) - A$
contained in  $E( K({\frac \alpha \beta}))$. If $E^\#  \cap D_{0} = \emptyset$ then since $E^\#  \cap 
\partial E(K) = \emptyset$ the disk  $E^\#$ is an inessential disk in this component  of $V_1 - A$ 
which is a solid torus.  If $E^\#  \cap D_{0} \neq \emptyset$ then since this solid torus is irreducible we can 
reduce the intersection  by isotoping $E^\#$  off the neighborhood of the half tunnel until $E^\#$  is
isotopic to $ D_{0}$.

\noindent b) If  $E$ is contained in the component $E( K({\frac \alpha \beta})) \cap V_1$
then, as above, since it is in the component of $V_1 - A$ which is a solid torus and cannot intersect
$A$ it is isotopic to $ D_{0}$ which is parallel to $D$.

\end{proof}

\vskip5pt

Assume in contradiction that the Heegaard splitting$(V_1, V_2)$ is weakly reducible and let  
${\cal D}_1,{\cal D}_2$ be a pair of essential disks in $V_1$ and $V_2$ respectively, so that
${\cal D}_1\cap {\cal D}_2 = \emptyset$. As $E( K({\frac \alpha \beta}))$ contains no interior
tunnel it follows from Corollary \ref{innertunnel} that all outermost disks of ${\cal D}_2 \cap A$  
are in $E(K_n)$. 

If the disk  ${\cal D}_1 \cap A = \emptyset $ then it is either contained in  
$E(K_n)$  or parallel to the disk  $D$: As if it is not in $E(K_n)$ it must be a non-essential 
disk in the solid torus $V_1^1$ and these are parallel to $D$.
In the first case it is essential in the strongly irreducible induced Heegaard splitting on $E(K_n)$ 
and so must intersect the outermost sub-disks of any essential disk ${\cal D}_2 \subset V_2$:  
Note that all outermost sub-disks of $ V_2$ which are contained in $E(K_n)$  are essential disks in 
the strongly irreducible Heegaard splitting induced on $E(K_n)$. In the second, case as all 
outermost sub-disks of $ V_2$ intersect the parallel copy of $D \subset E(K_n)$ it follows that
the correspomding disks of $V_2$ must run through the annulus $A$ and intersect $D = {\cal D}_1$.

If the disk ${\cal D}_1 \cap A \neq \emptyset $ then assume first, that the outermost sub-disk 
${\cal D}^\# \subset{\cal D}_1$  is in the $E(K_n)$ component  of $E(K) - A$.
By the above claim $\cal D^\#$ is an  essential disk there. Since the induced Heegaard splitting 
on$E(K_n)$  is strongly irreducible any two outermost sub-disks of ${\cal D}_1$ and ${\cal D}_2$ in $E(K_n)$
must intersect. 

If the outermost sub-disks of ${\cal D}_1$ are in $E(K({\frac \alpha \beta}))$ then by 
the claim above if we cut this component of $V_1$ along $\cal D^\#$ we obtain two components one 
of which is a solid torus and the other is a  $3$-ball $\cal B$ (see Fig. 8(a) and 8(b)).

\break

%\vglue-170pt

\vbox{\hskip30pt \epsfysize150pt\epsfbox{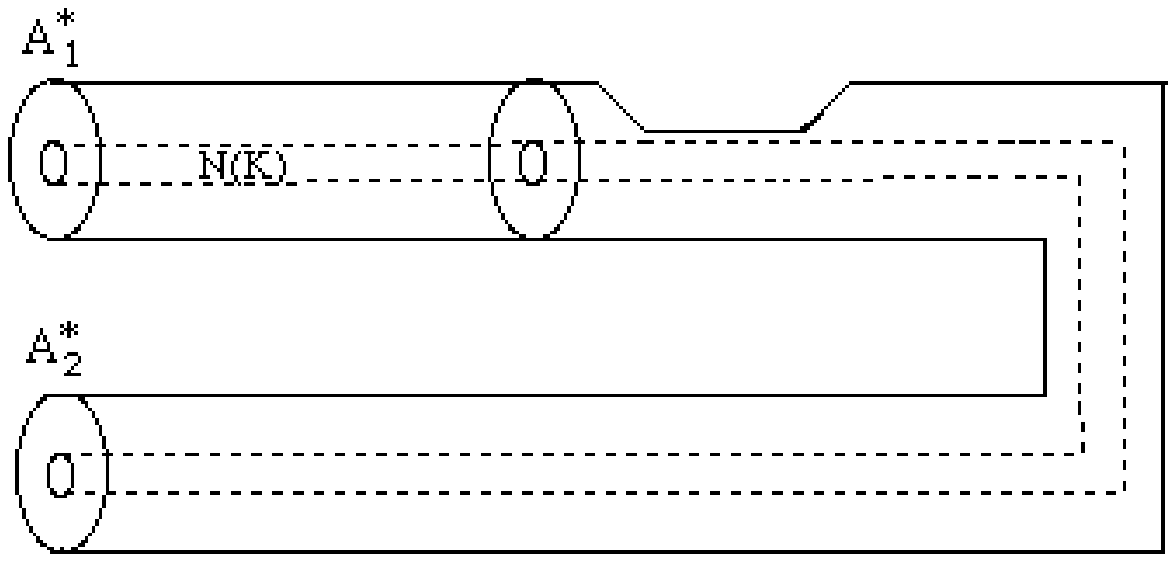}

\vskip15pt\centerline{Fig. 8(a)}}

\vskip10pt

\vbox{\hskip40pt \epsfysize180pt\epsfbox{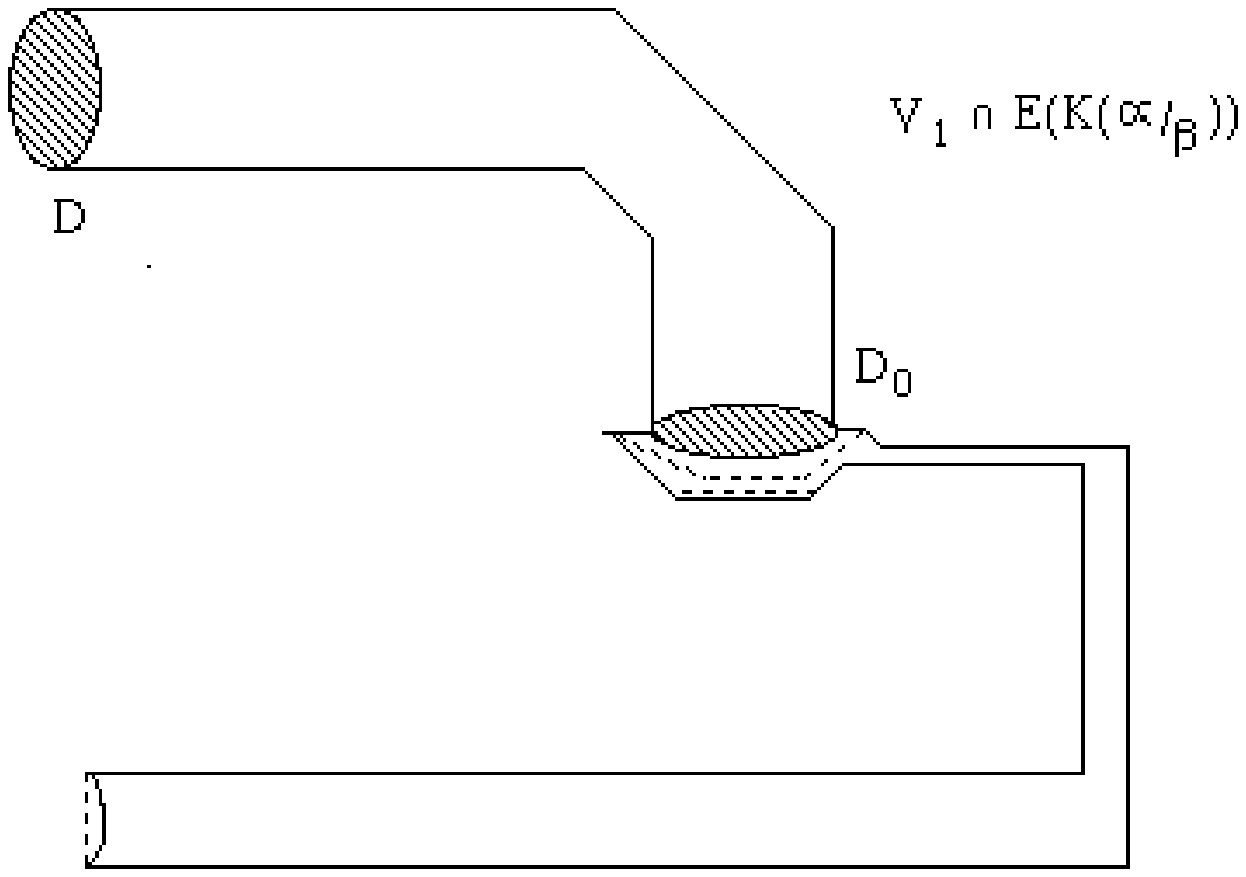}

\vskip15pt\centerline{Fig. 8(b)}}

\vskip20pt

Consider now an essential disk ${\cal D}_2$ in $V_2$.  If ${\cal D}_2 \cap A = \emptyset$ 
then ${\cal D}_2$ is an essential disk in $V_2^1$ or $V_2^2$, the two components of $V_2 - A$, depending 
on which side of $A$ the disk ${\cal D}_2$ is. Hence ${\cal D}_2$ is an essential disk in the
handlebody part of the induced Heegaard  splitting on either $E(K({\frac \alpha \beta}))$ or $E(K_n)$.
However these Heegaard splittings are strongly irreducible so ${\cal D}_2$ must intersect $D$ the cocore disk of
$t_1$ as it is an essential disk in the corresponding $V_1^1$ or  $V_1^2$. This implies that ${\cal D}_2$ must
intersect the decomposing annulus which is a contradiction. Hence ${\cal D}_2  \cap A $ is non-empty.

Let ${\cal D}^* \subset {\cal D}_2$ be a sub-disk, which is outermost among
all  sub-disks of ${\cal D}_2 - A$ which are contained in the  $E(K({\frac \alpha \beta}))$ 
component of $E(K) - A$. Let $\alpha_1, ..., \alpha_n$ be the components of 
${\cal D}^* \cap A$, then for all but one, say $\alpha_1$, the arcs $\alpha_i$ are ourtermost arcs of $D_2$ 
and hence are of type II (as in the proof of Theorem \ref{singtunthm}). Hence $\alpha_2, ..., \alpha_n$ have 
both end points on $D$, the cocore disk of the tunnel $t_1$. The arc $\alpha_1$ may be of type II 
or type I in which case  it has one end point on one of $\partial A_1^*$ or $\partial A_2^*$, and one on  
$\partial D$.
  
%\vglue-10pt

\vbox{\hskip40pt \epsfysize240pt\epsfbox{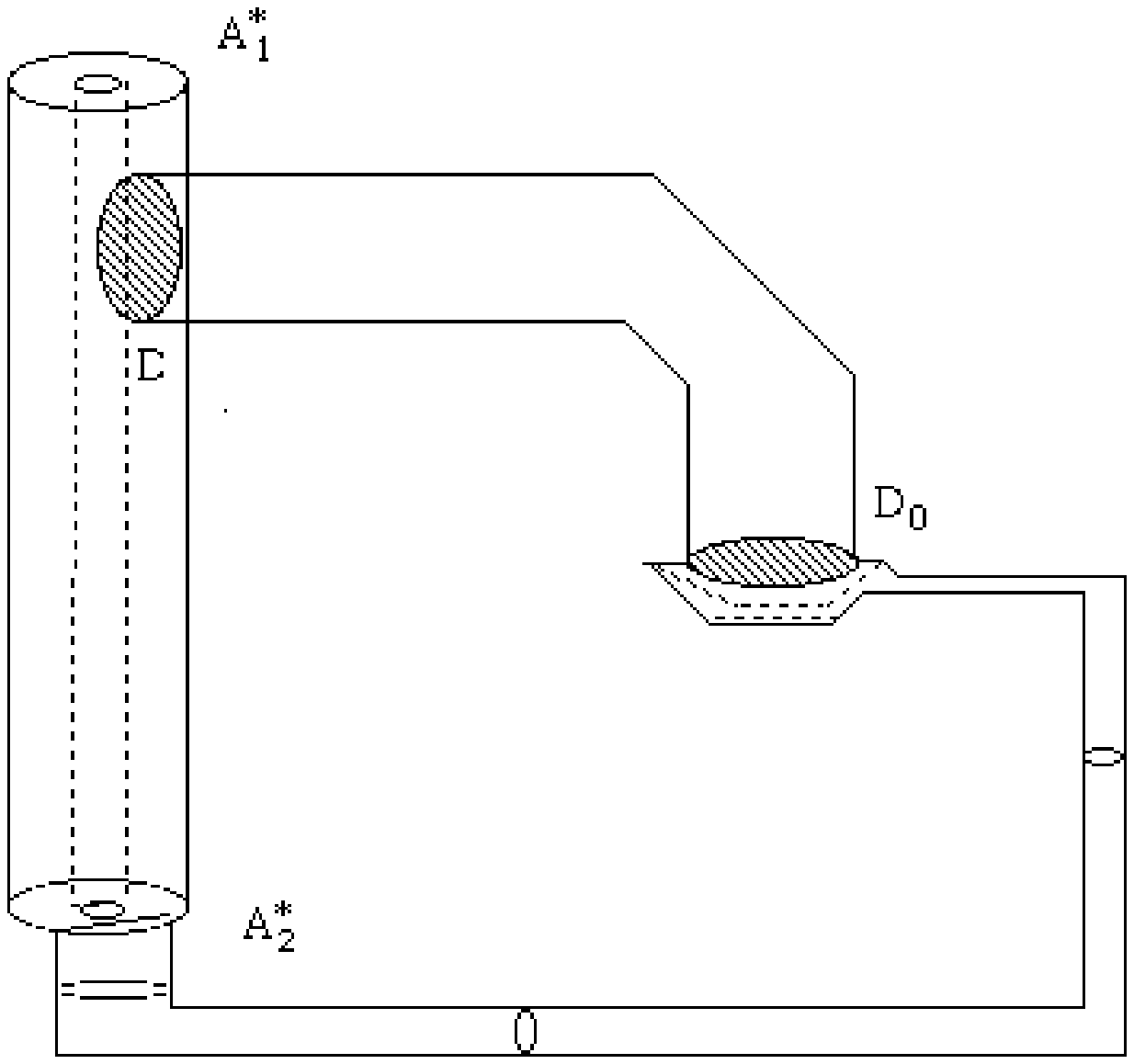}

\vskip15pt\centerline{Fig. 9}}

\vskip15pt

Since we are assuming that ${\cal D}_1 \cap {\cal D}_2 = \emptyset$, in both cases 
$\partial {\cal D}^* \cap F_1$ is a set of arcs contained in the annular sub-surface of $F_1$ 
depicted in Fig. 8b and Fig. 9 with all but at most one endpoint on  $\partial D$. Since 
by assumption all these arcs must be essential in $F_1$, it follows that $n = 1$ and  $\alpha_1$ 
is of type I. But this contradicts the fact that an outermost arc of intersection cannot be of type I
as then we can reduce the intersection of $A$ with $V_1$ in contradiction to the choice of $A$.
Thus  we have showed that any two essential disks in $V_1$ and $V_2$ must intersect and hence the 
Heegaard splitting $(V_1, V_2)$ is strongly irreducible.

\end{proof}

\vskip10pt

\begin{remark} \rm  The induced Heegaard splitting of genus three  on  $E(K_n \# K({\frac \alpha \beta}))$
is a stabilization of  the minimal Heegaard splitting $(V_1, V_2)$ discussed above. This can  be seen as follows: 
Remove a regular neighborhood of a short arc $\tau$ on $A$ connecting $\partial D$ to one of the vertical annuli, 
say $A_1^*$  from $V_2$ and add it as a $1$-handle to $V_1$. The arc $\tau$ is of type $I$ on some meridional 
disk $E$ of $V_2$ and since there is only one tunnel  crossing $A$ it bounds a sub-disk $\Delta$ on $E$. Hence the
cocore disk of $N(\tau)$ intersects $\Delta$ in a single point and therefore the pair $(V_1 \cup N(\tau), V_2 - N(\tau))$
is a stabilized Heegaard splitting for $E(K_n \# K({\frac \alpha \beta}))$. However we can slide the tunnel off $A$
by splitting it and sliding along $N(\tau)$. We obtain a isotopic Heegaard splitting with no tunnels crossing $A$
which is isotopic to the Heegaard splitting of $E(K_n \# K({\frac \alpha \beta}))$ which is induced by the two
 \lq\lq standard" Heegaard splittiings of $E(K_n)$ and $E( K({\frac \alpha \beta}))$.

\end {remark}

\vskip28pt

\section {References}
\vskip10pt

\noindent [CG] \hskip 20pt  A. Casson, C. Gordon; {\it Reducing Heegaard splittings}, Topology and 

\noindent \hskip 46pt  its Applications 27 (1987), 275 - 283.

\vskip8pt
\noindent [Go] \hskip 20pt  C. Gordon; {\it On primitive sets of loops in the boundary of a 

\hskip 28pt handlebody}, Topology and  its Applications 27 (1987), 285 - 299.

\vskip10pt

\noindent [Ha]\hskip 28pt  K. Hartshorn; {\it  Heegaard splittings of Haken manifolds have bounded

\noindent \hskip 48pt  distance}, preprint.

 \vskip10pt

\noindent [He]\hskip 28pt  J. Hempel; {\it 3-manifolds as viewed from the curve complex}, 

\noindent \hskip48pt Topology  40 (2001), 631 -- 657.

\vskip10pt  

\noindent [LM]\hskip 24pt  M. Lustig, Y. Moriah; {\it Closed incompressible surfaces in comp- 

\noindent \hskip 48pt lements of wide knots and links}, Topology and its Applications 92   

\noindent \hskip 48pt (1999), 1 - 13.   

\vskip10pt

\noindent [Me]\hskip 28pt  W. Menasco; {\it Closed incompressible surfaces in alternating knot

\noindent \hskip 48pt and link complements}, Topology 23 (1984), 37 - 44.

\vskip10pt 

\noindent [Mh]\hskip 27pt  Y. Moriah; {\it Incompressible surfaces and connected sum of knots}, J. 

\noindent \hskip 53pt  of Knot Theory and its Ramifications, 7 (1998), 955 - 965.

\vskip10pt 

\noindent [Mh1]\hskip 27pt  Y. Moriah; {\it On boundary primitive manifolds and a theorem of 

\noindent \hskip 53pt Casson-Gordon}, Topology and its Applications 125 (3) (2002), 

\noindent \hskip 53pt 571 -  579.

\vskip10pt

\noindent [Mo1]\hskip 23pt K. Morimoto; {\it  On the additivity of tunnel number of knots}, Topol- 

\noindent \hskip 53pt ogy and its Applications 53 (1993),  37 - 66.

\vskip10pt

\noindent [Mo2]\hskip 24pt K. Morimoto; {\it  There are knots whose tunnel numbers go down un- 

\noindent \hskip 51pt der connected sum}, Proc. Amer. Math. Soc. 123 (1995), 3527 - 

\noindent \hskip 52pt 3532.

\vskip10pt

\noindent [Mo3]\hskip 24pt K. Morimoto; {\it Tunnel number, connected sum and meridional

\noindent \hskip 51pt essential surfaces}, Topology 39 (2000), 469 - 487.

\vskip10pt

\noindent [Mo4]\hskip 24pt K. Morimoto; {\it On the super additivity of tunnel number of knots},

\noindent \hskip 53pt Math. Ann. 317 (2000), 489 - 508.

\vskip10pt

\noindent [MR]\hskip 28pt Y. Moriah, H.Rubinstein; {\it Heegaard Structures of negatively curved 

\noindent \hskip 53pt 3-manifolds}, Comm. in Anal. and Geom. 5 (1997), 375 - 412. 

\vskip10pt

\noindent [MS]\hskip 28pt K. Morimoto, J. Schultens; {\it Tunnel numbers of small knots do not 

\noindent \hskip 52pt  go down under connected sum},Proc. Amer. Math. Soc 128 (1) 

\noindent \hskip 52pt (2000), 269-278.

\vskip10pt

\noindent [MSY] \hskip 18pt K. Morimoto, M. Sakuma, Yokota; {\it Examples of tunnel number 

\noindent \hskip 53pt one knots which have the property "1 + 1 = 3" }, Math. Proc. 

\noindent \hskip 53pt Camb. Phil. Soc. 119 (1996), 113 - 118.

\vskip10pt

\noindent [Sc]  \hskip 28pt J. Schultens; {\it Additivity of tunnel number for small knots}, Preprint.

\vskip10pt

\noindent [Th]  \hskip 26pt A. Thompson; {\it The disjoint curve property and genus 2 manifolds},

\noindent \hskip 48pt Topology and its Applications 97 (1999), 273 - 279.

\vskip25pt

\obeyspaces   Yoav Moriah

\obeyspaces   Department of Mathematics

 \obeyspaces  Technion, Haifa  32000,

\obeyspaces  Israel

\vskip10pt

 ymoriah@tx.technion.ac.il

\end{document}